\documentclass[12pt]{article}
\usepackage[vmargin=1in,hmargin=1in]{geometry}

\usepackage{amssymb, cite}
\usepackage{graphicx}

\newcommand{\ee}{\end{equation}}
\newcommand{\bqn}{\begin{eqnarray}}

\newcommand{\eqn}{\end{eqnarray}}

\newcommand{\bd}{\begin{description}}
\newcommand{\ed}{\end{description}}
\newtheorem{Theorem}{Theorem}[section]
\newtheorem{lemma}[Theorem]{Lemma}

\newtheorem{claim}[Theorem]{Claim}

\newtheorem{stat}{}[section]

\def\bs{\begin{stat}}
\def\es{\end{stat}}
\def\ben{\begin{enumerate}}
\def\een{\end{enumerate}}

\def\bp{\noindent{\bf Proof}  \ }
\newcommand{\ep}{\hfill $\square$}

\makeatletter
\@addtoreset{equation}{section}

\makeatother

\begin{document}

\begin{center}
{\large {\bf GRAPHS WITH THE SAME
\\[2ex]
TRUNCATED CYCLE 
MATROIDS}}
\\[5ex]
 {\bf Jos\'e F. De Jes\'us  (University of Puerto Rico, San Juan,Puerto Rico, U.S.A.) }
\\[4ex]
 {\bf Alexander Kelmans (University of Puerto Rico, San Juan,Puerto Rico, U.S.A.) }

\end{center}

\date{}

\vskip 3ex

\begin{abstract}
The classical Whitney's $2$-Isomorphism Theorem describes the families of graphs having the same cycle matroid. 
In this paper we describe the families of graphs having the same truncated cycle matroid and prove, in particular, that every 3-connected graph, except for $K_4$, is uniquely defined by its truncated cycle matroid.

 \vskip 2ex


 \vskip 1ex


\end{abstract}

\section{Introduction}
\label{Intro}

\indent
  
 In 1933 H. Whitney \cite{W2} described the families of graphs having the same cycle matroid. He also proved, 
in particular, that every 3-connected graph is uniquely defined by its cycle matroid \cite{W3}.    
In this paper we describe the families of graphs having the same truncated cycle matroid and prove, in particular, that every 3-connected graph, except for $K_4$, is uniquely defined by its cycle matroid. The dual version of our paper is described by R. Chen and Z.Gao \cite{CG}.

\section{Main notions,  notation, and simple observations}
\label{Notations}

Given two finite sets $X$ and $Y$, let 
$X\Delta Y = (X\setminus Y) \cup (Y\setminus X)$.
 \\
 Given  a finite set $E$, let $2^E$ denote the set of all subsets of $E$ and 
  ${E\choose k}$ denote the set of all $k$-element subsets of $E$.
 \\
 Given  $ {\cal S} \subseteq 2^E$, let 
 $\Delta ({\cal S}) = \Delta \{S: S \in {\cal S}\}$.

\subsection{On graphs}
\label{graphs}
 
A {\em graph} $G$ is a triple $(V, E,\phi)$ such that $V$ and $E$ are disjoint finite sets, 
$V\cap E = \emptyset $, $V \ne \emptyset $, and 
$\phi: E\rightarrow {V\choose 2}$. We will also put 
$V = V(G)$ and $E = E(G)$.

The elements of $V = V(G)$ and $E = E(G) $ are called {\em vertices and  edges of graph} $G$, respectively. 
If $\phi(e) =  \{u,v\}$, we say that  
{\em vertices $u$ and $v$ are incident to edge $e$  and are the end-vertices of $e$ in $G$}.  

  \vskip 1.5ex
  
We say that {\em graph $G' = (V', E',\phi')$ is a subgraph of graph $G = (V, E,\phi)$ and write $G'  \le G$}
if $V'  \subseteq V$,  $E'  \subseteq E$, and function $\phi'$ is a restriction of  function $\phi'$ on $V'$. 

\vskip 1.5ex

The  {\em degree   $d(v,G)$ of vertex $v$ in $G$} is  the number of edges incident to $v$ in $G$. 

A {\em cycle in graph $G$} is an $ \le $-minimal subgraph of $G$ with every vertex of degree two. 
An {\em  $(x,y)$-path in graph $G$} is an $ \le $-minimal subgraph of $G$ with exactly  two vertices $x$ and $y$ of degree one.

A cycle $C$ in graph $G$ is called {\em Hamiltonian}, if $V(C) = V(G)$. A cycle $Q$ in graph $G$ is called 
{\em quasi-Hamiltonian}, if $V(Q) = V(G)\setminus q$ for some 
$q \in V(G)$.

A {\em forest in graph $G$} is a subgraph of $G$ with no cycles.
A forest $F$ in graph $G$ is called {\em maximal} if $F$ is not a proper subgraph of another forest in $G$. 

A {\em graph $G$ with at least  two vertices  is connected} if 
$G$ has an $(x,y)$-path for every two vertices $x$ and $y$ 
in $G$. Obviously, a maximal forest in a connected graph $G$ is 
a tree $T$ with $V(T) = V(G)$ (called a {\em spanning tree of $G$}). 

 \vskip 1.5ex

Let $G =(V, E,\phi)$ and $X \subseteq E$. 
Let $G[X]$ be the graph such that $E(G[X]) = X$ and $V(G[X])$ is the set of vertices of $G$ incident to at least one edge in $X$. We say that 
{\em $G[X]$ is the subgraph of $G$ induced by the edge subset 
$X$}.

   \vskip 1.5ex
 
 A graph $G$ is called {\em even} if every vertex of $G$ is incident to an even number of edges in $G$ (i.e. the degree of every vertex in $G$ is even). For example, a cycle is an even graph.
 
  \vskip 1.5ex

Graphs  $G =( V, E ,\phi)$ and $G' = (V', E', \phi')$ are {\em equal} if 
$V = V'$, 
$E = E'$, and $\phi = \phi'$.

\vskip 1.5ex

An {\em isomorphism from $G  = (V, E, \phi)$ to $G' = (V', E', \phi')$ } is a pair $(\nu, \varepsilon)$, where 
\\
$\nu: V \to V'$ and $\varepsilon: E \to E'$ are bijections such that
 $\phi(e)= \{x,y\} \Leftrightarrow \phi' (\varepsilon (e)) = \{ \nu(x),\nu(y)\}$.     
Graphs $G$ and $G'$ are {\em isomorphic} (denoted by $G \sim G'$) if there exists an isomorphism from $G$ to $G'$ (or, equivalently, an isomorphism from $G'$ to $G$). 

 \vskip 1.5ex
A {\em vertex star in graph $G$} is the set of edges in $G$ incident to the same vertex.
Graphs $G$ and $G'$ are {\em strongly isomorphic} (denoted by 
$G \approx G'$) if  they have the same family of vertex stars.

  \subsection{On matroids}
\label{matroids}

  Let $M$ be a matroid on the ground set $E$ with the set of bases 
 ${\cal B} = {\cal B}(M)$  and the set of circuits ${\cal C}(M)$. A circuit $H$ in $M$ is called {\em Hamiltonian} if the number of elements in $H$ is equal to the number of elements in a base of $M$ plus one.
\vskip 1.5ex
 We will need the following (known and easy to prove) fact.
 
 \begin{claim}
\label{Bt}
  Let $M$ be a matroid on the ground set $E$ and 
   
 ${\cal B}_t(M) = \{B - x:  B  \in  {\cal B}(M),  x  \in B \}$.
 
Then
\\
$(c1)$ 
${\cal B}_t$ is the set of bases of a matroid on the ground set $E$ {\em (denoted by $M_t$)} and
\\
$(c2)$ ${\cal C}(M_t) = {\cal C}(M)\cup {\cal B}(M) \setminus  {\cal H}(M)$,
where ${\cal H}(M)$ is the family of Hamiltonian circuits of $M$.
 \end{claim} 
 
 Matroid $M_t$ is called the the {\em truncation of matroid 
 $M$} 
 or, simply, a {\em trucated matroid}.

   \subsection{On matroids of a graph} 
\label{matroids,graphs}
 
 Let $G =( V, E ,\phi)$ be a graph with the vertex set $V$, the non-empty edge set $E$, and the incident function $\phi$. 
 Let ${\cal B}(G)$ be the family of the edge sets of  
 maximal forests and ${\cal C}(G)$ be the family of the edge sets
 of cycles in graph $G$.
 
 It is known and easy to prove that ${\cal B}(G)$ is the set of bases and ${\cal C}(G)$ is the set of circuits of 
 a {\em matroid $M$ on the ground set $E$} called {\em the cycle matroid of graph $G$} and denoted by $M(G)$. 
 
 Obviously, if graph $G$ is connected, then ${\cal B}(G)$ is 
 the family of the edge sets of spanning trees in $G$ 
 and ${\cal B}_t(G)$ is the family of the edge-sets 
 of maximal two-component forests in $G$.

 \vskip 1.5ex

\section{Main Result}  
\label{Main}

We   need the following (easy to prove) claim.
\begin{claim} 
\label{Delta} 
Suppose that
\vskip 1ex
\noindent
$(a1)$ 
$G$ is a graph and $Y, Z \subseteq E(G)$ and 
\vskip 1ex
\noindent
$(a2)$ 
$G[Y]$ is an even subgraph of $G$ and $G[Z]$ is a cycle of $G$.
\vskip 1ex
%\noindent
Then  
%\\
\vskip 1ex
\noindent
$(c1)$
$G[Y\Delta Z]$ is also an even subgraph of $G$, and therefore
\vskip 1ex
\noindent
$(c2)$ if ${\cal F}$ is a family of the edge sets of cycles in $G$, then 
$G[\Delta {\cal F}]$ is an even subgraph 
\vskip 1ex
\noindent
of $G$.
 \end{claim}

Put ${\cal C}(M_t(G)) ={\cal C}_t(G)$.
From the  definition of $M_t(G)$ we have:

\begin{lemma} 
\label{T,Q,S}
 Let $G$ be a connected graph with no loops and no parallel edges.
\vskip 1ex 
Then 
${\cal C}_t(G) = {\cal T}(G) \cup {\cal Q}(G) \cup {\cal S}(G)$, where
\vskip 1ex
\noindent
${\cal T}(G)$ is the family of the edge-sets of spanning trees of $G$,  
 \vskip 1ex
 \noindent
$ {\cal Q}(G)$ is the family of the edge-sets of quasi-Hamiltonian cycles of 
$G$, and 
 \vskip 1ex
 \noindent
$ {\cal S}(G)$  is the family of the edge-sets of small cycles of 
$G$.

\end{lemma}

We say that a pair of graphs $\{G, F\}$
satisfies  condition ${\cal K}$ if the following conditions hold:
\vskip 1ex
\noindent
$(b0)$ $G$ and $F$ are connected graphs,
 \vskip 1ex
 \noindent
$(b1)$ $G$ and $F$ are  graphs with no loops and no parallel edges,
  \vskip 1ex
\noindent
$(b2)$
$E(G) = E(F) = E$, 
 \vskip 1ex
 \noindent
$(b3)$
$M_t(G) = M_t(F)$, and
 \vskip 1ex 
 \noindent
$(b4)$
$M(G) \ne M(F)$.

\begin{lemma} 
\label{Lemma}
Suppose that a pair of graphs $\{G, F\}$
satisfies  condition ${\cal K}$.
 
Then there exists $X  \subseteq E$ such that  exactly one of the following holds:
\\
$(c1)$  $G[X]$ is a quasi-Hamiltonian cycle in $G$ and $F[X]$   is 
a spanning tree  in $F$ or
\\
$(c2)$  $F[X]$ is a quasi-Hamiltonian cycle in $F$ and $G[X]$ is 
 a spanning tree  in $G$.
 \end{lemma} 
   
\bp 
Put ${\cal C}(M_t(G)) ={\cal C}_t(G)$. 
By Lemma \ref{T,Q,S}, 
$
{\cal C}_t(G) 
= {\cal T}(G) \cup {\cal Q}(G) \cup {\cal S}(G)$.
\\
Obviously,
 \vskip 0.7ex
$M_t(G) = M_t(F)  \Leftrightarrow  {\cal C}_t(G) = {\cal C}_t(F)$. 
It is easy to see that 
$ {\cal S}(G)  = {\cal S}(F)$.  
\vskip 1ex 
 \noindent
Therefore, 
 \vskip 1ex
$M_t(G) = M_t(F)  \Leftrightarrow 
{\cal T}(G) \cup {\cal Q}(G) = {\cal T}(F) \cup {\cal Q}(F)$.
\vskip 1ex 
 \noindent
Now 
 \vskip 1ex
$M(G) \ne M(F) \Rightarrow {\cal T}(G) \ne {\cal T}(F) 
\Leftrightarrow {\cal T}(G) \Delta {\cal T}(F) \ne \emptyset $. 
\vskip 1ex 
 \noindent
Thus,
there exists $X \subseteq E$ such that at least one of the following holds:
\vskip 1ex 
 %\noindent
$X \in {\cal Q}(G)\setminus {\cal Q}(F)$
and 
$X  \in {\cal T}(F)\setminus {\cal T}(G)$ or 
\vskip 1ex 
 %\noindent
$X \in {\cal Q}(F)\setminus {\cal Q}(G)$
and 
$X  \in {\cal T}(G)\setminus {\cal T}(F)$.
\ep

 \begin{claim}
\label{XnotSymDf}
The edge set $X$  in Lemma \ref{Lemma} is not the symmetric difference of the edge sets of small cycles in $G$ $($ i.e. of members of set 
${\cal S}(G)$ $)$.
\end{claim}
\bp  
Suppose, on the contrary, that $X = \Delta {\cal F}$ for some family of edge sets of cycles in $G$.
Then by Claim $\ref{Delta}$, $F[\Delta {\cal F}]$ is an even graph and therefore is not a spanning tree in graph $F$. This contradicts 
Lemma $\ref{Lemma}$.
\ep

\vskip 3ex
Without loss of generality we can assume by Lemma  \ref{Lemma} that 
$G[X]$ is a quasi-Hamiltonian cycle in $G$ and $F[X]$   is 
a spanning tree  in $F$.

Let $G[X] = R$ and $F[X] = T$. Then $V(G) \setminus V(R)$ consists of one vertex that we denote by $c$. Since graph $G$ is connected, there exists 
an edge $s$ incident to $c$ and to a vertex (say, $r$) in $R$. 
Thus, $R' = R \cup (csr)$ is a connected spanning subgraph of  
$G$, and so 
if $e$ is an edge in $G$, then one of  its end-vertices is in $R$ and the other one is either equal to $c$ or is also in $R$.  

We call cycle $R$ the {\em rim} of $G$, an edge of  $G$ incident to 
$c$  a {\em spoke} of $G$, and an edge $h$ in $E(G) \setminus E(R)$  with both end-vertices in $R$ a {\em chord} of $G$.

 \begin{claim}
\label{chord}
Graph $G$ has no chords, and therefore $G$ has at least one spoke.
\end{claim} 

\bp  
 %\vskip 0.7ex
 Suppose, on the contrary, that $G$ has a chord $h$ with the end-vertices  
 $u$ and $v$. Then $h$ belongs to  two cycles $C_1$ and  $C_2$ in $G$ such that 
\vskip 1ex 
 \noindent
 $(p1)$
 $E(C_1), E(C_2) \in {\cal S}(G) = {\cal S}(F)$ and
\vskip 1ex 
 \noindent
$(p2)$
 $E(C_1)\Delta E(C_2) = E(R)$.
 \vskip 1ex 
% \noindent
This contradicts Claim \ref{XnotSymDf}.
\ep

 \vskip 3ex

%Let $s_1, \ldots, s_n$ be the list of spokes in $G$ and $r_i$ a vertex in $R$ incident to   spoke $s_i$. 

\begin{claim}
\label{some2spokes}
Suppose that $s$ and $s'$ are some different spokes in $G$ with the end-vertices $r$ and $r'$ in $R$, respectively. Let $P = P(r,r')$ be a shortest path in $R$ with the end-vertices $r$ and $r'$. Then path $P$ has at most two edges.
\end{claim}

\bp   
Let $H = rscs'r'$ and  let $P' = P'(r,r')$ be the path in $R$ such that 
$P \cup P' = R$.
 Let $C = P \cup H$ and $C' = P' \cup H$. Then $E(R) = E(C) )\Delta E(C')$. 
Suppose, on the contrary, that $p = e(P) \ge 3$. 
Then also $e(P')  \ge 3$, and so $e(R) \ge 2p \ge 6$.
Also $e(C) = e(P) + 2$ and $e(C') = e(P') + 2$. 
Therefore $e(C) < e(R)$ and $e(C') < e(R)$. 
Then,  $E(C), E(C') \in  {\cal S}(G)$ and $E(C)\Delta E(C') = E(R)$.
Now by Lemma $\ref{Delta}$, $F[E(C)\Delta E(C')]$ is an even graph and therefore is not a spanning tree in graph $F$. This contradicts 
Lemma $\ref{Lemma}$.
\ep

\vskip 1ex 
\noindent
\begin{claim}
\label{<4spokes}
Graph $G$ has at least one and at most three spokes. Moreover,
\vskip 1ex
\noindent
 $(c1)$ if $G$ has two spokes $s_1$ and $s_2$, then their end-vertices $r_1$ and $r_2$ in $R$ are  on distance one or two in $R$
(see Fig. 3.7[c1,1] and Fig.3.7[c1,2]) 
  and
\vskip 1ex
\noindent
$(c2)$ if $G$ has three spokes  $s_1$, $s_2$, and $s_3$, then their 
end-vertices $r_1$, $r_2$, and $r_3$ belong to a 3-vertex path in $R$
(see Fig3.7[c2,p1]).
\end{claim}

\bp  
By Claim \ref{chord}, $G$ has at least one spoke.
If $G$ has at most two spokes, then our claim $(c1)$ follows from  
Claim \ref{some2spokes}. 
%\\

Now we prove $(c2)$.
Suppose  that  $G$ has exactly three spokes
$s_1$, $s_2$, and $s_3$ with the end-vertices $r_1$, $r_2$, and 
$r$ in $R$, respectively.
\vskip 1ex 
\noindent
$(p1)$
Suppose first that some two vertices from $\{r_1, r_2, r\}$, say, $r_1$ and $r_2$, are adjacent in $R$. If $r$ is adjacent to $r_1$ or $r_2$, then 
$r_1$, $r_2$, and $r$ belong to a 3-vertex path in $R$, and
we are done (see Fig.[c2,p1]). 
Therefore $r$ is on distance at least 2 from $r_1$ and $r_2$ in $R$.
Then $E(R) = X$ is the symmetric difference of the edge sets of three small cycles. This contradicts Claim \ref{XnotSymDf}.
\vskip 1ex 
\noindent
$(p2)$
Now suppose that no two vertices from $\{r_1, r_2, r\}$ are adjacent in $R$.  Then by Claim \ref{some2spokes}
every two vertices from $r_1$, $r_2$, $r$ are on distance 2 in $R$. Therefore $R$ is a 6-edge cycle and $E(R) = X$ is the symmetric difference of 
the edge sets of three squares  $S_i$, $i \in \{1,2,3\}$ in $G$.
This contradicts Claim  \ref{XnotSymDf}.
\ep

\vskip 1.5ex
Put $M_t = M_t(G) = M_t(F)$. 
 
\begin{claim}
\label{Isomorphism1}
Suppose that 
\vskip 1ex
\noindent
 $(a1)$ a pair of graphs $\{G, F\}$ satisfies  condition ${\cal K}$ and
 \vskip 1ex
\noindent
 $(a2)$ $G$ has exactly one spoke $s$ with its end-vertex $r$ in $R = G[X]$. 
\vskip 1ex
Then 
\vskip 1ex
\noindent
 $(c1)$
 either $F$ is isomorphic to $G$ or $F$ is a cycle and 
\vskip 1ex
\noindent
 $(c2)$ 
$M_t $ is isomorphic to a uniform matroid 
%\vskip 1ex
%\noindent
$U_ {m-2,m} = (E, {\cal B})$,
\vskip 1ex
\noindent
where $m$ is the number of elements in $E$ and 
${\cal B}$ = $ {E}\choose{m-2}$   
is the set of bases of matroid $U_ {m-2,m}$.
\end{claim}

\bp  
As we have assumed above, $R = G[X]$ is a quasi-hamiltonian cycle in $G$ and $T= F(X)$ is a spanning tree in $F$.
Then $T\cup s = F[X \cup s]$ has a cycle, say $D$. If $D$ is a small cycle, then $E(D)$ induces a small cycle in $G$, a contradiction. Therefore one of the following holds:  
\vskip 1ex
\noindent
 $(h1)$ 
 $D$ is a quasi-hamiltonian cycle in $F$, 
 and so $F$ is isomorphic to $G$ or 
\vskip 1ex
\noindent
 $(h2)$ 
 $D$ is a hamiltonian cycle in $F$, and so $D = F$ is a cycle.  
\vskip 1ex
Therefore $(c1)$ holds.
 Moreover, in both cases $M_t$ satisfies $(c2)$.
 \ep

\begin{claim}
\label{Isomorphism2,1}
Suppose that
\vskip 1ex
\noindent
 $(a1)$
a pair of  graphs $\{G, F\}$ satisfies  condition ${\cal K}$ and
\vskip 1ex
\noindent
 $(a2)$
$G$ has exactly two spokes $s_1$ and $s_2$ and their 
 end-vertices $r_1$ and $r_2$ are incident to an edge, 
 say $x$, in 
 $R = G[X]$. 
\vskip 0.7ex
Then 
\vskip 1ex
\noindent
 $(c1)$
graphs $G$ and $F$ are isomorphic,
\vskip 1ex
\noindent
 $(c2)$ graph $F$ can be obtained from graph $G$ by a series of some Whitney 2-vertex cut switches and, possibly, also by exchanging two edges $e$ and $s_i$ for some $i \in \{1,2\}$ (which is not a Whitney graph operation),
 and
\vskip 1ex
\noindent
 $(c3)$
 graphs $G$ and $F$ are 2-connected.

\end{claim}

\bp  
By Claim \ref{Lemma},
 $R = G[X]$  is a quasi-hamiltonian cycle in $G$ and $F[X]$ is a spanning tree $T$ in $F$.  
Since $G$ has exactly two spokes $s_1$ and $s_2$ and their 
end-vertices $r_1$ and $r_2$ in $R = G[X]$ are  on distance one in 
$R$,
there is   an edge $x$ in $R$  such that 
$\triangle  =  cs_1r_1x r_2s_2c$  is a  3-cycle in $G$. 
Clearly, $e \in E(R) \cap E(T)$.
Then $E(\triangle)$ is the edge-set of a  3-cycle $Z$ in $F$.

Let $V(Z) = \{z, z_1, z_2 \}$, where vertices  $z_1$, $z_2$ are
incident to edge $x$,  and so vertex $z$ is not incident to edge $x$.
Let $zPz'$ be a shortest path in tree $T$ from $z$ to $\{z_1, z_2 \}$.  Let 
$z' = z_1$,  and so $s_1$ is the edge in $Z$ incident to vertices $z$ and  $z_1$ in $T$. Obviously, $x \not \in E(P)$. If $e(P) < e(T - x)$, then $P \cup s_1$  is 
a small cycle in $F$, and therefore $E(P) \cup s_1$ is the edge set  of a small cycle in $G$ distinct from  $\triangle $, 
a contradiction. Then  
$e(P) = e(T - x)$, and so $E(P) = E(T - x)$, 
a tree $T$ is the path $zPz_1xz_2$, and $F = T  \cup s_2$.
Thus, $F$ is isomorphic to $G$, and therefore $(c1)$ holds. 
It also follows that $(c2)$ and $(c3)$ holds.
\ep

\begin{claim}
\label{Isomorphism2,2}
Suppose that 
\vskip 1ex
\noindent
 $(a1)$ a pair of graphs $\{G, F\}$ satisfy  condition ${\cal K}$ and
\vskip 1ex
\noindent
 $(a2)$
$G$ has exactly two spokes $s_1$ and $s_2$ and their 
 end-vertices $r_1$ and $r_2$ in $R$ are  the end-vertices 
 of a  two-edge path $r_1x_1rx_2r_2$ in $R$.
 
\vskip 1ex
Then 
\vskip 1ex
\noindent
 $(c1)$
graphs $G$ and $F$ are isomorphic,
\vskip 1ex
\noindent
 $(c2)$ graph $F$ can be obtained from graph $G$ by 
 a series of some Whitney 2-vertex cut switches and, possibly, 
 also by exchanging two edges $e_i$ and $s_j$ for some 
 $i, j \in \{1,2\}$ 
 (which is not a Whitney graph operation), 
 and
\vskip 1ex
\noindent
 $(c3)$
 graphs $G$ and $F$ are 2-connected.
\end{claim}

\bp  
By $(a2)$, 
$\diamondsuit  =  
cs_1r_1x_1rx_2r_2s_2c$  is an induced 4-cycle in $G$. 
Then $E(\diamondsuit)$ is also the edge-set of an induced  
4-cycle in $F$.
Let $R^*$ be the cycle in $G$ with 
$E(R^*) = (E(R) \setminus \{x_1,x_2\}) \cup \{s_1,s_2\}$.

Consider the path  $Z = R - r$ with the end-vertices $r_1$ and $r_2$.
Then  $Z = R^* - c$.
Obviously, $e(Z) \ge 1$.  Also  
$\{s_1,s_2,x_1,x_2\}$ is the edge set of a 4-cycle in $F$.

Suppose first that $Z$ has one edge, say $z$. Then 
$Q_c = (R  \cup z)  \setminus r$ and $Q_r = R^* \cup z)  \setminus c$ are triangles and the only quasi-hamiltonian cycles in $G$.
By Lemma \ref{Lemma}, we can assume that  $E(Q_c)$ induces a spanning tree in $F$. Then 
$v(F) = v(G) = 4$, $E(F) = E(G)$, and $e(F) = e(G) = 5$. 
Therefore graphs  $G$ and $F$ are isomorphic.

Now suppose that  $e(Z) \ge 2$.
Let $z_1, z_2 \in E(Z)$,  where $z_1 \ne z_2$. 
Let  
$Z_i = Z - z_i$ and $T_i = Z_i \cup (r_1s_1cs_2r_2)\cup \{x_i, r\}$, 
where $i \in \{1,2\}$. 
Then each $T_i$ is a spanning tree of $G$ and 
$G = T_1  \cup T_2$. 
Let $E_i = E(T_i)$.
Obviously, $F = F(E_1) \cup F(E_2)$.

We claim that each $F(E_i)$ is also a spanning tree of $F$.
Suppose, on the contrary, that, say,   $F(E_1)$ is not a spanning tree of $F$. Then $D = F(E_1)$ is a quasi-hamiltonian cycle  in $F$  containing  a path with three edges $x_1$, $s_1$, $s_2$. Hence $x_2$ is a cord of cycle $D$.
Then  the edge-set 
$(E(D)\cup x_2)\setminus \{s_1,s_2,x_1\}$ induces
 a small cycle in $F$ but not in $G$, a contradiction.
It follows that graphs $G$ and $F$ are isomorphic. 
Therefore $(c1)$ holds.
It follows that $(c2)$ and $(c3)$  also hold.
 \ep

\begin{claim}
\label{Isomorphism3} 
Suppose that
\vskip 1ex
\noindent
$(a1)$ a pair of graphs $\{G, F\}$ satisfies  condition ${\cal K}$ and 
\vskip 1ex
\noindent
 $(a2)$ $G$ has exactly three spokes  $s_1$, $s_2$, and $s$ 
 such that end-vertices $r_1$, $r_2$, and $r$ 
 belong to a 3-vertex path
$r_1x_1rx_2r_2$ in $R$ and spoke $s$ is incident to $r$.

\vskip 1ex
Then 
\vskip 1ex
\noindent
 $(c1)$
graphs $G$ and $F$ are isomorphic,
\vskip 1ex
\noindent
 $(c2)$ graph $F$ can be obtained from graph $G$ by a series of some Whitney 2-vertex cut switches and, possibly, also by exchanging two edges $e_i$ and $s_i$ for some $i \in \{1,2\}$
 (which is not a Whitney graph operation), 
 and
\vskip 1ex
\noindent
 $(c3)$
 graphs $G$ and $F$ are 2-connected.
\end{claim}

\bp 
 From  $(a2)$ it follows (as above) that  in  graph $G$
there is a vertex $r$ and edges $x_1$ and $x_2$ in $R$  such that 
$\diamondsuit  =  cs_1r_1x_1rx_2s_2r_2c$  is a 4-cycle in $G$.
Let $G' = G \setminus s$ and $F' = F \setminus s$. 
Then by Claim \ref{Isomorphism2,2}, $G'$ and $F'$ are isomorphic. 
If $F = F' \cup s$ is not isomorphic to $G$, then 
$s$ belongs to a small cycle $D$ in $F$, 
which is also a small cycle in $G$, a contradiction. 
Therefore  $(c1)$ holds. It also follows that  $(c2)$ holds.
\ep

\vskip 1ex 

It is easy to prove the following  
\begin{claim}
\label{HnotConnected}
Let $H$ be a non-connected graph with no loops, no parallel edges, no isolated vertices, and 
 with the set of components $\{C_i: i = 1, ..., k\}$. 
 Let $V = \{v_i \in V(C_i): i = 1, ..., k\}$ and
 let $H*$ be the  graph obtained from $H$ by identifying the vertices in $V$ with a new vertex $v$. 

 Then 
 \vskip 1ex
\noindent
 $(c1)$ $v(H^*) \ge 3$ and $H^*$ is of  connectivity one
 $($i.e. is connected, but not 2-connected$)$,
\vskip 1ex
\noindent
 $(c2)$$M(H) = M(H^*)$, 
\vskip 1ex
\noindent
 $(c3$ $M_t(H) = M_t(H^*)$,and
 \vskip 1ex
\noindent
 $(c4$ $H^*$ satisfies condition ${\cal K}$.
 
\end{claim}

\vskip 1ex 
\noindent

%If $A$ is a connected graph, we will also put sometimes 
% $A= A^*$.

\begin{claim}
\label{G,FnotConnected} 
Suppose that a pair of graphs $\{G, F\}$ satisfies all  condition in ${\cal K}$, 
except  for condition  $(b0)$, 
namely, at least  one of  $G$, $F$, say  G, 
is not a connected  graph.
Then $G$ has two components, namely, a one-edge component and 
a cycle.
\end{claim}

\bp 
By Claim \ref{HnotConnected}, $G^*$ satisfies condition 
${\cal K}$. 

By Claims \ref{Isomorphism1}, \ref{Isomorphism2,1},    \ref{Isomorphism2,2}, and \ref{Isomorphism3}, the graphs satisfying the conditions of these Claims are 2-connected. 
Since graph $G^*$ is not 2-connected, 
$G^*$ does not satisfy the assumptions of these Claims. 
Thus, $G^*$ satisfies the assumptions of Claim \ref{Isomorphism1}. Therefore $G$ has exactly two components, namely, a one-edge component and a cycle.
\ep

From the above we have in particular:
\begin{claim}
\label{Gis3-Connected} 
Every 3-connected graph, except for $K_4$, is uniquely defined by its truncated cycle matroid.
\end{claim}

\end{document}